\documentclass[10pt]{amsart}

\usepackage{amsmath}
\usepackage{amsthm}
\usepackage{hyperref}

\newtheorem{prop}{Proposition}
\newtheorem{lemma}[prop]{Lemma}
\newtheorem{thm}[prop]{Theorem}
\newtheorem{cor}[prop]{Corollary}
\newtheorem{conj}[prop]{Conjecture}
\newtheorem*{mainthm}{Main Theorem}

\theoremstyle{definition}
\newtheorem{defn}[prop]{Definition}

\newcommand{\dt}{\frac{\partial}{\partial t}}
\newcommand{\ddt}{\frac{d}{d t}}
\newcommand{\brs}[1]{\left| #1 \right|}
\newcommand{\gD}{\Delta}

\newcommand{\gS}{\Sigma}
\newcommand{\gl}{\lambda}
\newcommand{\gt}{\theta}
\renewcommand{\ge}{\epsilon}
\newcommand{\N}{\nabla}
\newcommand{\FF}{\mathcal F}
\renewcommand{\bar}[1]{\overline{#1}}

\newcommand{\hook}{\mathbin{\hbox{\vrule height2.4pt width4.5pt depth-2pt
\vrule height5pt width0.4pt depth-2pt}}}

\DeclareMathOperator{\Rc}{Rc}
\DeclareMathOperator{\Vol}{Vol}

\begin{document}

\title{Ricci Yang-Mills flow on surfaces}

\author{Jeffrey Streets}

\address{Fine Hall\\
         Princeton University\\
         Princeton, NJ 08544}
\email{\href{mailto:jstreets@math.princeton.edu}{jstreets@math.princeton.edu}}

\thanks{Author was partially supported by a Clay Liftoff Fellowship and by the
National Science Foundation via DMS-0703660}

\date{July 31st, 2009}

\begin{abstract}  We study the behaviour of the Ricci
Yang-Mills flow for $U(1)$
bundles on surfaces.  We show that existence for the flow reduces to a bound on
the isoperimetric constant.  In the presence of such a bound, we show that on
$S^2$, if the
bundle is nontrivial, the
flow exists for all time.  For higher genus surfaces the flow always exists for
all time.  The volume normalized flow always exists for all time and converges
to a constant scalar curvature metric with the bundle curvature $F$ parallel. 
Finally, in an appendix we classify all gradient solitons of this flow on
surfaces.
\end{abstract}



\maketitle

\section{Introduction}
Fix $(M^n, g)$ a Riemannian manifold.  Suppose $L \to M$ is the
total space of a $U(1)$-bundle over $M$, and $A$ is a connection on
this bundle with curvature $F$.  This $F$ is a purely imaginary two-form on $M$
which represents the first Chern class of the line bundle associated to $L$.  In
what follows we will often not refer to
the total space of the bundle and focus attention on $M, g$ and $A$, and
furthermore identify $F$ with a real valued two-form.  We
say that a family $(M, g(t), A(t))$ is a solution to Ricci
Yang-Mills flow (RYM-flow) if
\begin{gather}
\begin{split}
\dt g_{ij} =&\ -2 \Rc_{ij} + g^{kl} F_{ik} F_{jl}\\
\dt A =&\ - d^* F.
\end{split}
\end{gather}
This equation was studied in \cite{thesis} in the hope that by
introducing connections $A$ where the curvature $F$ has special
properties then the flow would have behaviour simpler than that of the Ricci
flow.  Moreover, this system
arises naturally in
physics as the renormalization group flow for a certain nonlinear sigma model. 
Also, a recent paper of Lebrun \cite{Lebrun} shows an interesting connection
between solutions to the static equation, known as the Einstein-Maxwell
equation, to the existence of extremal K\"ahler metrics in dimension $4$.
Finally, we mention that this equation has generated interest as a tool for
better understanding magnetic flows on surfaces \cite{Dan}.  By examining
homogeneous solutions, the
following conjecture is plausible:

\begin{conj} Let $(M^{2n}, g)$ be a Riemannian manifold.  Let $L \to M$
be the total space of a $U(1)$ bundle over $M$.  Given $A$ a
connection on $L$ satisfying
\begin{gather*}
[F^{\wedge n}] \neq 0
\end{gather*}
then the solution to the Ricci Yang-Mills flow with initial
condition $(g, A)$ exists for all time.
\end{conj}
\noindent We mention a related conjecture for odd-dimensional manifolds in the
conclusion. In this paper we examine this conjecture in the
case $n =
1$.  We show that the regularity of the flow can be reduced to showing a bound
on the Sobolev constant of the manifold.  Recall that the Sobolev constant of a
Riemannian surface $(M^2, g)$ is the smallest constant $C_S$ such that the
inequality
\begin{gather} \label{Sobolev}
\left(\int_M \brs{f - \bar{f}}^2 dV_g\right)^{\frac{1}{2}} \leq C_S \int_M
\brs{\N f}
\end{gather}
holds for any function $f \in C^1(M)$, where $\bar{f}$ is the average value of
$f$.  It is known that this constant is equivalent to other Sobolev constants,
and moreover is equivalent to the isoperimetric ratio \cite{Chavel}.

\begin{thm} \label{existence} Let $g$ be a Riemannian metric on an oriented
surface $M$, and let $L \to M$ denote the total space of a $U(1)$-bundle over
$M$ with connection $A$.  Let $(g(t), A(t))$ be the solution to RYM flow with
this condition.  If the solution goes singular at time $T < \infty$, then either
\begin{align*}
\lim_{t \to T} \Vol(g(t)) = 0
\end{align*}
or
\begin{align*}
\lim_{t \to T} C_S(g(t)) = \infty.
\end{align*}
If $(g(t), A(t))$ is the solution to volume-normalized RYM flow and it goes
singular at time $T < \infty$, then
\begin{align*}
\lim_{t \to T} C_S(g(t)) = \infty.
\end{align*}
In other words, in the presence of volume and Sobolev constant bounds the
solution to RYM flow on a surface is nonsingular. 
\end{thm}

We of course expect the isoperimetric constant to stay bounded along the flow. 
In fact, such a bound for the Ricci flow on $S^2$ was shown by Hamilton
\cite{Hamilton}.  Such a bound for solutions to RYM flow is as yet unclear.
 We are moreover able to completely describe the limiting behaviour of
infinite-time solutions with no hypotheses.  The overall situation is described
in the main theorem below.

\begin{mainthm} \label{maintheorem} Suppose solutions to RYM flow satisfy a
Sobolev constant bound.  In other words, given $(M^{2}, g, A)$ a solution to RYM
flow on $[0, T]$, one has $C_S(g(t)) < C(T)$ for all $t \leq T$.  Let $g$ be a
Riemannian metric on
an oriented surface $M$, and let $L \to M$ denote the total space of a
$U(1)$-bundle over $M$ with connection $A$.
\begin{enumerate}
\item{ If $M \cong S^2$ and $[F] \neq 0$, the solution to RYM flow with initial
condition $(g, A)$ exists for all time.  If the volume stays finite at infinity,
the solution converges to the round metric with $F$ parallel.  Moreover, the
volume normalized flow exists for all time and converges to the round metric
with
$F$ parallel.}
\item{ If $M \cong S^2$ and $[F] = 0$, the solution to volume normalized RYM
flow with initial
condition $(g, A)$ exists for all time and converges to the round metric with $F
\equiv 0$.}
\item{ If $\chi(M) \leq 0$, the solution to RYM flow with initial condition $(g,
A)$ exists for all time and the volume-normalized flow exists for all time and
converges to a constant curvature metric with $F$
parallel.}
\end{enumerate}
\end{mainthm}
In fact, all of the convergence statements for flows existing for all time hold
without the Sobolev constant bound hypothesis.  Notice that two important
questions are left unresolved.  In particular, we do
not know if the unnormalized equation on $S^2$ with $[F] \neq 0$ has a volume
bound and hence converges at infinity.  Also, it would be interesting to know
the complete behaviour on $S^2$ when $[F] = 0$.  We conjecture that the solution
goes singular in finite time, converging to a round point with $F \equiv 0$. 
Since understanding gradient solitons may play a role in resolving these issues,
we provide a classification.
\begin{prop} If $g$ is a gradient soliton on a closed
surface $\gS^2$ then $g$ has constant curvature and $F$ is parallel.
\end{prop}

The proof of long time existence in the presence of the Sobolev constant bound
generalizes the corresponding proof for Ricci
flow found by Struwe
\cite{struwe}.  We first reduce to a flow on a conformal factor $u$ and a
connection $A$, and indeed we show that a certain energy functional generalizing
the
Liouville energy for the conformal factor to include the Yang-Mills coupling is
monotonically
decreasing along a solution to RYM flow.  Using this and a further a priori
integral estimate we are able to bound the $H^2$ norms of $u$ and $A$, and thus
prove long time existence.  The Moser-Trudinger inequality plays a key role in
the proof as well. Given the long time existence, we are able to show that the
Calabi energy remains bounded, and thus apply the compactness result of Xiuxiong
Chen \cite{Chen} to show convergence at infinity.  We note that Andrea Young has
independently obtained stability results for the Ricci Yang-Mills flow on a
surface \cite{AY}.

In section 2 we rewrite the RYM flow equation on a surface in terms of a
conformal factor and introduce the volume normalized equation.  Section 3
contains certain
a-priori integral estimates, section 4
completes the proof of Theorem \ref{existence}, and section 5 has the proofs of
convergence, completing the proof of the Main Theorem.  Section 6 is a
concluding discussion, and section 7
is an appendix containing the classification of Ricci Yang-Mills
solitons on surfaces.

The author would like to thank Gang Tian for helpful discussions and for
suggesting the approach of \cite{struwe}.  Thanks also go to Dan Jane for
several stimulating conversations.  Finally, the author would like to thank an
anonymous referee for a very thorough reading of and many valuable comments on
an earlier version of this manuscript.

\section{Reduction to Conformal Flow}
In this section we show that the metric component of the Ricci Yang-Mills flow
on a surface is a
conformal flow.
We already know that on a surface $\Rc = \frac{1}{2} R g$.  On a
surface the term $g^{kl} F_{ik} F_{jl}$ is a scalar multiple of the metric as
well.
\begin{lemma} Given $(M^2, g)$ a Riemannian surface and $F \in
\bigwedge^2 T^* M$, $g^{kl} F_{ik} F_{jl} = \frac{1}{2} \brs{F}_g^2
g_{ij}$.
\begin{proof} Fix a point $x \in M$ and choose normal coordinates for $g$ at
$x$.
In these coordinates we have $F(x) = \lambda(x) dx^1 \wedge dx^2$.
Clearly then
\begin{gather*}
g^{kl} F_{ik} F_{jl}(x) = \left(
\begin{matrix}
\lambda^2 & 0\\
0 & \lambda^2
\end{matrix} \right)
= \lambda^2(x) g(x).
\end{gather*}
Since the left and right hand sides are both tensors it follows that
there exists a function $\gl(x)$ so that $g^{kl} F_{ik} F_{jl} =
\gl^2(x) g_{ij}$.  Taking the trace of this equation gives $\gl^2(x) =
\frac{1}{2} \brs{F}_g^2$.
\end{proof}
\end{lemma}
\noindent Using this lemma the RYM-flow on a surface becomes the system of
equations
\begin{gather} \label{RYMSflow}
\begin{split}
\dt g =&\ - R g + \frac{1}{2} \brs{F}_g^2 g\\
\dt A =&\ - d^* F
\end{split}
\end{gather}
Furthermore, if $g = e^u g_0$ where $g_0$ is a fixed metric of constant
curvature $R_0$ and unit volume, then we can write
\begin{gather*}
R = e^{-u} \left(R_0 - \gD u \right)
\end{gather*}
where $\gD$ is with respect to the metric $g_0$.  Thus we can
write the RYM-flow as the system
\begin{gather} \label{confflow}
\begin{split}
\dt u =&\ e^{-u} \left( \gD u - R_0 + \frac{1}{2} e^{-u} \brs{F}^2
\right)\\
\dt A =&\ - d^* F
\end{split}
\end{gather}
where in this equation the norm $\brs{F}^2$ is taken with respect to
$g_0$.  Note that since $A$ is a connection on a $U(1)$ bundle $F$ may be
thought of as just a usual (closed) $2$-form on $M^2$.  Thus we derive the
evolution equation
\begin{gather} \label{Fflow}
\begin{split}
\dt F =&\ \dt \left(d A \right)\\
=&\ - d d^*_{g} F\\
=&\ \gD_{d, g} F\\
=&\ \gD_g F
\end{split}
\end{gather}
where $\gD_{d, g}$ is the Laplace-Beltrami operator of $g$.  Note that the
curvature term in the B\"ochner formula for $n$-forms always vanishes on
$n$-manifolds \cite{wu}, thus the last line follows where $\gD_g$ is the rough
Laplacian of $g$.  We take the time here to mention an important convention in
this paper.  Any metric which is used without further decoration will be the
fixed background metric.  Any time we use the time-dependent metric $g(t)$ we
will decorate the quantity with a $g$.

We will also need a certain volume-normalized system.  Note that
\begin{align*}
\int_M R_g dV_g - \frac{1}{2} \int_M \brs{F}^2_g dV_g = R_0 -
\frac{1}{2} \int_M e^{-u} \brs{F}^2 dV
\end{align*}
Thus consider
\begin{gather} \label{vnflow}
\begin{split}
\dt u =&\ e^{-u} \gD u + R_0(1 - e^{-u}) + \frac{1}{2} \left(e^{-2u}
\brs{F}^2 - \frac{\int_M e^{-u} \brs{F}^2 dV}{\int_M e^u dV} \right)\\
\dt A =&\ - d^* F
\end{split}
\end{gather}
The volume of $g(t)$ remains constant under this
evolution equation.  Note that this system does \emph{not} differ from the
unnormalized equation by a rescaling in space and time.  This is a consequence
of the fact that $\dt g$ does not have homogeneous scaling.  In particular, the
term $\frac{1}{2} \brs{F}^2_g g$ has inverse scaling with respct to the metric
while $R_g g$ has neutral scaling with the metric.  Also, $F$ still obeys
(\ref{Fflow}) with respect to this time dependent metric.

\section{Integral Estimates}
In this section we will prove a-priori integral estimates for the RYM-flow on a
surface.  First we define a
functional which is monotonic for solutions to RYM-flow
\begin{gather} \label{Fdef}
\FF(u,A) := \int_M \left(\brs{d u}^2 + e^{-u} \brs{F}^2 \right) dV + 2 R_0
\int_M u dV
\end{gather}
where the norms and volume form are those of the background metric $g_0$.
\begin{prop} \label{surfaceenergy} Given $(M^2, u(t), A(t))$ a solution to
(\ref{confflow}) we have
\begin{gather} \label{energycalc}
\begin{split}
\ddt \FF(u(t), A(t)) =&\ - 2 \int_M e^u \brs{u_t}^2 dV -2 \int_M \brs{\N^g
F}_g^2
dV_g
\end{split}
\end{gather}
\begin{proof} First we compute
\begin{align*}
\ddt \int_M \brs{d u}^2 dV =&\ 2 \int_M \left< d \left( e^{-u} \left(
\gD u - R_0 + \frac{1}{2} e^{-u} \brs{F}^2 \right) \right), d u
\right> dV\\
=&\ - \int_M \left(2 e^{-u} (\gD u)^2 - 2 R_0 e^{-u} \gD u + e^{-2u}
\brs{F}^2 \gD u \right) dV\\
=&\ 2 R_0 \int_M e^{-u} \brs{d u}^2 dV - 2 \int_M e^{-u} (\gD u)^2 dV
- \int_M e^{-2u} \brs{F}^2 \gD u dV
\end{align*}
Next we use the equation $\int_M e^{-u} \brs{F}^2 dV = \int_M
\brs{F}^2_{g} dV_{g}$ and compute using (\ref{RYMSflow}) and (\ref{Fflow})
\begin{align*}
\ddt \int_M \brs{F}^2_g dV_g =&\ -2 \int_M \brs{ \N^g F}^2_g dV_g +
\int_M \left( R_g - \frac{1}{2} \brs{F}_g^2 \right) \brs{F}^2_g
dV_g\\
=&\ -2 \int_M \brs{ \N^g F}^2_g dV_g - \int_M e^{-2u} \left( \gD u -
R_0 \right) \brs{F}^2 dV\\
&\ - \frac{1}{2} \int_M e^{-3 u} \brs{F}^4 dV
\end{align*}
Next we have
\begin{align*}
2 R_0 \ddt \int_M u dV =&\ 2 R_0 \int_M e^{-u} \gD u - R_0 e^{-u} + \frac{1}{2}
e^{-2 u} \brs{F}^2 dV\\
=&\ 2 R_0 \int_M e^{-u} \brs{d u}^2 - R_0 e^{-u} + \frac{1}{2} e^{-2 u}
\brs{F}^2 dV.
\end{align*}
Combining these calculations gives
\begin{align*}
\ddt &\FF (u(t), A(t))\\
=&\ 4 R_0 \int_M e^{-u} \brs{d u}^2 dV - 2 \int_M
\brs{ \N^g F}^2_g dV_g - 2 R_0^2 \int_M e^{-u} dV\\
&\ + 2 R_0 \int_M e^{-2u} \brs{F}^2 d V - 2 \int_M e^{-u} (\gD u)^2 dV\\
&\ - 2 \int_M e^{-2u} \brs{F}^2 \gD u dV -
\frac{1}{2} \int_M e^{-3 u} \brs{F}^4 dV\\
=&\ -2 \int_M e^u \brs{ e^{-u} \gD u - e^{-u} R_0 + \frac{1}{2} e^{-2 u}
\brs{F}^2}^2 dV - 2 \int_M \brs{\N^g F}^2_g dV_g\\
=&\ - 2 \int_M e^u \brs{u_t}^2 dV - 2 \int_M \brs{\N^g F}_g^2 dV_g
\end{align*}
as required.
\end{proof}
\end{prop}

\begin{prop} \label{volnormenergy} Given $(M^2, u(t), A(t))$ a solution to
(\ref{vnflow}) we have
\begin{gather} \label{volnormenergyev}
\begin{split}
\ddt \FF(u(t), A(t)) =&\ - 2 \int_M e^u \brs{u_t}^2 dV - 2 \int_M \brs{\N^g
F}^2_g dV_g
\end{split}
\end{gather}
\begin{proof} Adding a constant to $u$ clearly does not affect the evolution of
$\int_M \brs{d u}^2 dV$.  Next in computing the evolution of $\int_M e^{-u}
\brs{F}^2 dV$ we pick up 
\begin{gather*}
\left(\frac{1}{2} \int_M e^{-u} \brs{F}^2 dV  - R_0 \right) \int_M e^{-u}
\brs{F}^2 dV.
\end{gather*}
Thus from the previous proposition we compute
\begin{gather*}
\begin{split}
\ddt \FF(u(t), F(t)) =&\ -2 \int_M e^u \brs{ e^{-u} \gD u - e^{-u} R_0 +
\frac{1}{2} e^{-2 u} \brs{F}^2}^2 dV\\
&\ - 2 \int_M \brs{\N^g F}^2_g dV_g + 2 R_0^2 - 2 R_0 \int_M e^{-u} \brs{F}^2
dV\\
&\ + \frac{1}{2} \left( \int_M e^{-u} \brs{F}^2 dV \right)^2\\
=&\ - 2 \int_M e^u \brs{u_t}^2 - 2 \int_M \brs{\N^g F}^2_g dV_g
\end{split}
\end{gather*}
\end{proof}
\end{prop}

\begin{cor} \label{Fbound} Given $(M^2, u(t), A(t))$ a solution to
(\ref{confflow}) or (\ref{vnflow}) we have
\begin{gather*}
\FF(u(t), A(t)) \leq \FF(u(0), A(0))
\end{gather*}
\begin{proof} This follows from the above lemmas.
\end{proof}
\end{cor}

\begin{lemma} \label{volumelemma1} Let $(S^2, g(t), A(t))$ be a solution to the
Ricci
Yang-Mills flow on $S^2$ satisfying $[F] \neq 0$.  Then there exists
a constant $C = C(g_0, \brs{[F]})$ so that the inequality
\begin{gather}
\Vol(g(t)) \geq C
\end{gather}
holds for all time that the flow exists.
\begin{proof}  First note that on a Riemannian surface $(M, g)$, any $F \in
\bigwedge^2 T^* M$
satisfies $F = \pm \brs{F} dV$.  This implies the inequality
\begin{align*}
0 <&\ \brs{[F]}\\
=&\ \brs{\int_M F}\\
\leq&\ \int_M \brs{F}_g dV_g\\
\leq&\ \left( \int_M \brs{F}_g^2 dV_g \right)^{\frac{1}{2}}
\Vol(g)^{\frac{1}{2}}.
\end{align*}
Using this and the Gauss-Bonnet Theorem we compute the evolution
equation
\begin{align*}
\ddt \Vol(g(t)) =&\ - \int_M R dV + \frac{1}{2} \int_M \brs{F}^2 dV\\
\geq&\ - 4 \pi + \frac{\brs{[F]}^2}{2 \Vol(g(t))}.
\end{align*}
If $\Vol(g(t)) \leq \frac{[F]^2}{8 \pi}$ then $\dt
\Vol(g(t)) \geq 0$ and the result follows.
\end{proof}
\end{lemma}

\begin{lemma} \label{volumelemma2} Given $(M^2, g(t), A(t))$ a solution to
(\ref{confflow}) there exists a constant $C > 0$ depending on $(g(0), A(0))$ so
that the inequality
\begin{gather} \label{volumebound}
\Vol(g(t)) \leq Vol(g(0)) + C t
\end{gather}
holds for any $t > 0$.
\begin{proof} Using Proposition \ref{surfaceenergy} we estimate
\begin{align*}
\ddt \Vol(g(t)) =&\ \ddt \int_M dV_g\\
=&\ \int_M \left( -R + \frac{1}{2} \brs{F}_g^2 \right) dV_g\\
=&\ - 2 \pi \chi(M) + \frac{1}{2} \int_M \brs{F}_g^2 dV_g.
\end{align*}
However, using the fact that $\mathcal F$ is bounded and the Liouville energy is
bounded below in any conformal class we see
\begin{align*}
\int_M \brs{F}_g^2 dV_g =&\ \mathcal F(u(t), A(t)) - \int_M \left(\brs{d u}^2 +
2 R_0 u \right) dV\\
\leq&\ \mathcal F(u(0), A(0)) + C.
\end{align*}
Therefore
\begin{align*}
\ddt \Vol(g(t)) \leq C
\end{align*}
and the result follows.
\end{proof}
\end{lemma}

\begin{lemma} \label{energybound} Given $(M^2, g(t), A(t))$ a solution to
(\ref{confflow}) or (\ref{vnflow}), on any finite time interval $[0, T]$ there
exists a constant $C$ depending only
on $(u(0), A(0))$ and $T$ so that 
\begin{gather*}
\brs{\brs{\N u}}_{L^2} \leq C, \qquad \int_M e^{-u} \brs{F}^2 dV \leq C
\end{gather*}
\begin{proof} In the case $R_0 \leq 0$ using Jensen's inequality and the volume
bound we easily conclude that
\begin{gather*}
\int_M u \leq \log \int_M \left( e^u dV \right) \leq C.
\end{gather*}
Thus we have an a-priori lower-bound for $\FF$ and the result follows.  For the
case $R_0 > 0$ we must modify our flow by an explicit M\"obius transformation as
in \cite{struwe}.  Specifically we solve for $\phi(t)$ a family of conformal
diffeomorphisms of the sphere such that $h(t) = \phi(t)^* g(t) = e^{v(t)} g_0$
satisfies
\begin{gather*}
\int_M x dV_h = 0
\end{gather*}
for all time, where $x$ is the position vector in $\mathbb R^3$.  Note that
these diffeomorphisms are certainly different from those obtained for fixing the
conformal gauge of Ricci flow.  Since $\FF$ is
diffeomorphism invariant it follows from Proposition \ref{volnormenergy} (in the
volume normalized case) that $\FF$ is uniformly bounded for the
diffeomorphism-modified flow and thus in particular
\begin{gather*}
\int_M \brs{dv}^2 dV + 2 R_0 \int_M v dV < C.
\end{gather*}
Then using Aubin's result \cite{Aubin} and the volume bound we conclude
\begin{gather*}
\brs{\brs{v}}_{H^1}^2 \leq C
\end{gather*}
One now easily gets an $C^1$ bound on the diffeomorphism parameter $\phi$ as in
\cite{struwe} Lemma 6.2 which gives the requisite bounds.
\end{proof}
\end{lemma}

\begin{lemma} \label{MosTrud} Given $(M^2, g(t), A(t))$ a solution to
(\ref{confflow}) or (\ref{vnflow}), on any finite time interval $[0, T]$ for any
$k$
we have
\begin{gather}
\sup_{0 \leq t < T} \int_M e^{k \brs{u}} dV < \infty
\end{gather}
\begin{proof} Since the volume is bounded on any finite time interval by Lemma
\ref{volumelemma2} and $\brs{\brs{\N u}}_{L^2}$ is bounded on a finite time
interval by Lemma \ref{energybound}, the result follows from the Moser-Trudinger
inequality .
\end{proof}
\end{lemma}

\section{Long Time Existence}
In this section we will use the integral estimates of the previous section and
the assumed bound on the Sobolev constant to
get an $H^2$ bound for both $u$ and $A$.  These bounds prove Theorem
\ref{existence} and the existence statements of the Main Theorem.  In the next
section we will use the gradient property to get the
convergence statements of the Main Theorem.  We point out that a
general short-time existence theorem for RYM
flow was shown in \cite{thesis} using the DeTurck gauge fixing procedure for
both the Ricci flow and the Yang Mills flow together.  Our bounds will apply to
any flow whose volume is bounded over any
finite time interval.  In particular these estimates work to show long time
existence for the
volume normalized flow, and the unnormalized flow in the cases when $\chi(M)
\leq 0$ and when $\chi(M) > 0, [F] \neq 0$ by Lemma \ref{volumelemma1}.  We will
explicitly work with the
unnormalized flow.

We will make use of the multiplicative Sobolev inequality
\begin{gather} \label{multsob}
\brs{\brs{f}}_{L^4}^2 \leq C \brs{\brs{f}}_{L^2} \brs{\brs{f}}_{H^1} \leq C
\brs{\brs{f}}_{H^1}^2.
\end{gather}
The constant of this inequality is equivalent to the Sobolev constant as we have
defined it using H\"older's inequality.  Also we use an inequality of
Calder\'on-Zygmund type:
\begin{gather} \label{caldzyg}
\int_M \brs{\N^2 f}^2 dV \leq C \int_M \brs{\gD f}^2 dV.
\end{gather}
We will have occasion to write certain terms using the metric $g$ for notational
convenience, and we will mostly apply the Sobolev inequality with respect to the
fixed background metric.  There is one term which requires the use of the
Sobolev inequality for $g$, and we treat it explicitly.  Also, we will make
repeated use of Lemmas
\ref{energybound} and \ref{MosTrud}.

We start with a preliminary observation.  Since $\FF(u(t), A(t))$ is continuous,
nonincreasing and bounded below, given $\ge > 0$ there is a $\tau > 0$ so that
given any $0 \leq t_0 < t_1 \leq T$ such that $t_1 - t_0 < \tau$ we have
\begin{gather} \label{intest0}
\FF(u(t_0), A(t_0)) - \FF(u(t_1), A(t_1) \leq \ge
\end{gather}
In particular for such times one has the estimate
\begin{gather} \label{intest1}
\int_{t_0}^{t_1} \int_M \brs{\N^g F}_g^2 dV_g \leq \int_{t_0}^{t_1} \dt
\FF(u(t), A(t)) \leq \ge
\end{gather}
which follows from Proposition \ref{surfaceenergy}.  Consider the calculation
\begin{align*}
\ddt \int_M e^u \brs{\N^g F}_g^2 dV_g =&\ 2 \int_M e^u \left< \N_i^g \gD_g F,
\N^g_i
F \right>_g dV_g - \int_M \left( e^u \right)_t \brs{\N^g F}^2_g dV_g\\
=&\ 2 \int_M e^u \left< \N^g_j \N^g_j \N^g_i F + \left(R_g + F^{*2} \right) *
\N^g F + \N u_t * F , \N_i^g F \right>_g dV_g\\
&\ - \int_M \left( e^u \right)_t \brs{\N^g F}^2_g dV_g\\
=&\ - 2 \int_M e^u \brs{\N^g \N^g F}_g^2 dV_g + \int_M e^u \N u * \N^g \N^g F *
\N^g F\\
&\ + \int_M e^u \left(R_g + F^{*2} \right)* \N^g F * \N^g F dV_g + \int_M e^u \N
u_t * F * \N^g F\\
&\ - \int_M \left( e^u \right)_t \brs{\N^g F}_g^2 dV_g\\ 
\leq&\ - \int_M e^u \brs{\N^g \N^g F}^2_g dV_g\\
&\ + C \int_M \left(e^u \brs{\N u}_g^2 + \gD u + \brs{R_0} + e^{u} \brs{F}_g^2
\right) \brs{\N^g F}_g^2 dV_g\\
&\ + C \int_M e^u \brs{\N u_t}_g \brs{F}_g \brs{\N^g F}_g dV_g.
\end{align*}
In the second line we commuted derivatives and in the third line integrated by
parts.  First we estimate
\begin{gather} \label{intest105}
C \brs{R_0} \int_{t_0}^{t_1} \int_M \brs{\N^g F}_g^2 dV_g dt \leq C
\end{gather}
by (\ref{intest1}).  Now we estimate
\begin{align*}
\int_M e^u \brs{\N u}_g^2 \brs{\N^g F}_g^2 dV_g =&\ \int_M e^u \brs{\N u}^2
\brs{\N^g F}^2_g dV\\
\leq&\ \brs{\brs{\N u}}_{L^4}^2 \brs{\brs{ e^{\frac{u}{2}} \brs{\N^g
F}_g}}_{L^4}^2\\
\leq&\ C \brs{\brs{\N u}}_{L^2} \brs{\brs{\N u}}_{H^1} \brs{\brs{e^{\frac{u}{2}}
\brs{\N^g F}_g}}_{L^2} \brs{\brs{e^{\frac{u}{2}} \brs{\N^g F}_g}}_{H^1}\\
\leq&\ C \sup_{t_1 \leq t < t_2} \brs{\brs{u}}_{H^2}\brs{\brs{e^{\frac{u}{2}}
\brs{\N^g F}_g}}_{L^2} \brs{\brs{e^{\frac{u}{2}} \brs{\N^g F}_g}}_{H^1}
\end{align*}
Which implies
\begin{gather} \label{intest101}
\begin{split}
\int_{t_0}^{t_1} \int_M &e^u \brs{\N u}_g^2 \brs{\N^g F}_g^2 dV_g\\
\leq&\ C \sup_{t_0 \leq t < t_1} \brs{\brs{u}}_{H^2} \left( \int_{t_0}^{t_1}
\int_M \brs{\N^g F}_g^2 dV_g \right)^{\frac{1}{2}} \cdot \left( \int_{t_0}^{t_1}
\int_M \brs{\N e^{\frac{u}{2}} \brs{\N^g
F}_g}^2 dV \right)^{\frac{1}{2}}\\
\leq&\ C \ge \sup_{t_0 \leq t < t_1} \brs{\brs{u}}_{H^2}^2 + C \ge
\int_{t_0}^{t_1} \int_M e^u \left(\brs{\N^g \N^g F}_g^2 + \brs{\N u}^2_g
\brs{\N^g F}^2_g \right) dV_g\\
\leq&\ C \ge \sup_{t_0 \leq t < t_1} \brs{\brs{u}}_{H^2}^2 + C \ge
\int_{t_0}^{t_1} \int_M e^u \brs{\N^g \N^g F}_g^2 dV_g.
\end{split}
\end{gather}
Analogously to the above estimate we get
\begin{align*}
\int_M e^{u} \brs{F}_g^2 \brs{\N^g F}_g^2 dV_g \leq&\ \left(\int_M e^{2 u}
\brs{F}_g^4 dV \right)^{\frac{1}{2}} \left( \int_M e^{2 u} \brs{\N^g F}^4_g dV
\right)^{\frac{1}{2}}\\
\leq&\ C \brs{\brs{e^{\frac{u}{2}}\brs{F}_g}}_{L^2} \brs{\brs{ e^{\frac{u}{2}}
\brs{F}_g}}_{H^1} \brs{\brs{e^{\frac{u}{2}} \brs{\N^g F}_g}}_{L^2}
\brs{\brs{e^{\frac{u}{2}} \brs{\N^g F}_g}}_{H^1}\\
\leq&\ C \brs{\brs{ e^{\frac{u}{2}} \brs{F}_g}}_{H^1} \brs{\brs{e^{\frac{u}{2}}
\brs{\N^g F}_g}}_{L^2} \brs{\brs{e^{\frac{u}{2}} \brs{\N^g F}_g}}_{H^1}\\
\end{align*}
Integrating this in time and arguing as in (\ref{intest101}) yields
\begin{gather} \label{intest110}
\begin{split}
\int_{t_0}^{t_1} \int_M &\ e^u \brs{F}^2_g \brs{\N^g F}^2_g dV_g\\
&\ \leq C \ge \sup_{t_0 \leq t < t_1} \int_M e^u \brs{\N^g F}^2_g dV_g + C \ge
\int_{t_0}^{t_1} \int_M e^{u} \brs{\N^g \N^g F}_g^2 dV_g.
\end{split}
\end{gather}
Also we have the estimate
\begin{align*}
\int_{t_0}^{t_1} \int_M \brs{\N^g F}^2_g \gD u dV_g \leq&\ \int_{t_0}^{t_1}
\brs{\brs{\gD u}}_{L^2} \brs{\brs{e^{\frac{u}{2}} \brs{\N^g F}_g}}_{L^4}^2\\
\leq&\ C \ge \sup_{t_0 \leq t < t_1} \brs{\brs{u}}_{H^2}^2 + C \ge
\int_{t_0}^{t_1} \int_M e^u \brs{\N^g \N^g F}_g^2 dV_g.
\end{align*}
Turning to the final term, we see
\begin{align*}
\int_{t_0}^{t_1} \int_M & e^u \brs{\N u_t}_g \brs{F}_g \brs{\N^g F}_g dV_g dt\\
=&\ \int_{t_0}^{t_1} \int_M e^{-u} \brs{\N u_t} \brs{F} \brs{\N^g F} dV\\
\leq&\ \int_{t_0}^{t_1} \left(\int_M e^u \brs{\N u_t}^2 dV dt
\right)^{\frac{1}{2}} \left( \int_M e^{-3 u} \brs{F}^2 \brs{\N^g F}^2 dV
\right)^{\frac{1}{2}} dt\\
\leq&\ C \ge \int_{t_0}^{t_1} \int_M e^u \brs{\N u_t}^2 dV dt + \int_{t_0}^{t_1}
\int_M e^{u} \brs{F}_g^2 \brs{\N^g F}_g^2 dV_g dt\\
\leq&\ C \ge \int_{t_0}^{t_1} \int_M e^u \brs{\N u_t}^2 dV dt + C \ge
\int_{t_0}^{t_1} \int_M e^{u} \brs{\N^g \N^g F}_g^2 dV_g dt.
\end{align*}
where in the last line we applied (\ref{intest110}).  Combining these estimates
gives
\begin{gather} \label{intest2}
\begin{split}
\int_{t_0}^{t_1} &\int_M e^u \brs{\N^g \N^g F}^2_g dV_g + \sup_{t_0 \leq t <
t_1}
\int_M e^u \brs{\N^g F}_g^2 dV_g\\
\leq&\ C \ge \sup_{t_0 \leq t < t_1} \brs{\brs{u}}_{H^2}^2 + C \int_M e^u
\brs{\N^g F}_g^2 dV_g (t_0) + C
\end{split}
\end{gather}

Now we turn to estimating $u$.  Our bounds here are directly adopted from
section 6 of \cite{struwe}.  First we have
\begin{gather*}
\dt e^u - \gD u = - R_0 + \frac{1}{2} e^{-u} \brs{F}^2
\end{gather*}
Multiplying this equation by $- \gD u_t$ and integrating gives
\begin{align*}
\int_M &e^ u \brs{\N u_t}^2 dV + \frac{1}{2} \dt \int_M \brs{\gD u}^2 dV\\
\leq&\ \frac{1}{2} \int_M e^{u} \brs{\N u_t}^2 dV + C \int_M e^u \brs{\N u}^2
\brs{u_t}^2 dV -
\int_M e^{-u} \brs{F}^2 \gD u_t dV
\end{align*}
Integrating in time and using the estimate
\begin{gather*}
\brs{\brs{u}}_{L^2}^2 \leq C \int_M e^{2 \brs{u}} dV \leq C(T),
\end{gather*}
which follows from Jensen's inequality and Lemma \ref{MosTrud}, we conclude
\begin{gather} \label{intest3}
\begin{split}
I :=&\ \int_{t_0}^{t_1} \int_M e^{u} \brs{\N u_t}^2 dV dt + \sup_{t_0 \leq t <
t_1} \brs{\brs{u}}_{H^2}^2\\
&\ \leq C \int_{t_0}^{t_1} \int_M e^u \brs{\N u}^2 \left( \brs{u_t}^2 + 1
\right) dV dt\\
&\ - \int_{t_0}^{t_1} \int_M e^{-u} \brs{F}^2 \gD u_t dV dt +
\brs{\brs{u(t_0)}}_{H^2} + C.
\end{split}
\end{gather}
Since $e^u$ is bounded in $L^2$ we deduce from the Sobolev inequality
\begin{gather} \label{intest5}
\int_M e^{u} \brs{\N u}^2 dV \leq C \brs{\brs{\N u}}_{L^4}^2 \leq C(T)
\brs{\brs{u}}_{H^2}^2 \leq C(T) \sup_{t_0 \leq t < t_1} \brs{\brs{u}}_{H^2}^2
\end{gather}
Similarly, using the Sobolev inequality (\ref{multsob}) and the a-priori bound
on $\brs{\brs{\N u}}_{L^2}$ we are able to bound
\begin{gather} \label{intest10}
\begin{split}
\int_M e^u \brs{\N u}^2 \brs{u_t}^2 dV =&\ 4 \int_M \brs{\N u}^2 \brs{\left(
e^{\frac{u}{2}} \right)_t}^2 dV\\
\leq&\ C \brs{\brs{ \N u}}_{L^4}^2 \brs{\brs{\left( e^{\frac{u}{2}}
\right)_t}}_{L^4}^2\\
\leq&\ C \brs{\brs{\N u}}_{L^2} \brs{\brs{\N u}}_{H^1} \brs{\brs{ \left(
e^{\frac{u}{2}} \right)_t}}_{L^2} \brs{\brs{ \left( e^{\frac{u}{2}}
\right)_t}}_{H^1}\\
\leq&\ C \brs{\brs{ \left( e^{\frac{u}{2}} \right)_t}}_{L^2} \brs{\brs{ \left(
e^{\frac{u}{2}} \right)_t}}_{H^1} \sup_{t_0 \leq t < t_1} \brs{\brs{u}}_{H^2}
\end{split}
\end{gather}
We need to estimate the time integral of the first two terms in the above
expression.  First of all it is clear that
\begin{gather} \label{intest20}
\int_{t_0}^{t_1} \brs{\brs{ \left( e^{\frac{u}{2}} \right)_t}}_{H^1}^2 dt \leq C
\int_{t_0}^{t_1} \int_M e^u \left( \brs{\N u_t}^2 + \brs{\N u}^2 \brs{u_t}^2 +
\brs{u_t}^2 \right) dV dt
\end{gather}
To estimate the other integral, we use (\ref{energycalc}) to compute
\begin{gather*}
\begin{split}
\brs{\brs{ \left( e^{\frac{u}{2}} \right)_t}}_{L^2}^2 =&\ 
\int_M e^u \brs{u_t}^2 dV\\
=&\ - \frac{d}{d t} \FF(u(t), A(t)) - \int_M \brs{ \N^g F}^2_g dV_g\\
\leq&\ - \frac{d}{d t} \FF(u(t), A(t)).
\end{split}
\end{gather*}
Thus we can conclude
\begin{gather} \label{intest30}
\begin{split}
\int_{t_0}^{t_1} \brs{\brs{ \left( e^{\frac{u}{2}} \right)_t}}_{L^2}^2 dt \leq&\
\FF(u(t_0), A(t_0)) - \FF(u(t_1), A(t_1))
\end{split}
\end{gather}
Thus integrating (\ref{intest10}) in time, applying H\"older's inequality and
using (\ref{intest20}) and (\ref{intest30}) gives
\begin{gather*}
\begin{split}
II :=& \int_{t_0}^{t_1} \int_M e^u \brs{\N u}^2 \brs{u_t}^2 dV dt\\
\leq&\ C \left( \FF(u(t_0), A(t_0)) - \FF(u(t_1), A(t_1)) + t_1 - t_0
\right)^{1/2}(I + II + C)
\end{split}
\end{gather*}
Now, recall from (\ref{intest0}) that we can choose $t_1 - t_0$ small enough
that
\begin{gather*}
C \left(\FF(u(t_0), A(t_0)) - \FF(u(t_1), A(t_1)) + t_1 - t_0 \right)^{1/2} \leq
\ge \leq \frac{1}{2}
\end{gather*}
which implies
\begin{gather*}
II \leq 2 \ge I + C
\end{gather*}
Thus from (\ref{intest3}) and (\ref{intest5}) we conclude
\begin{gather} \label{intest40}
I \leq C \left(t_1 - t_0 + \ge \right) I + C \brs{\brs{u(t_0)}}_{H^2}^2 -
\int_{t_0}^{t_1} \int_M e^{-u} \brs{F}^2 \gD u_t dV dt + C(T)
\end{gather}

We now turn to the last term in this expression.  First of all by integration by
parts and the Cauchy-Schwarz inequality we have
\begin{gather} \label{intest50}
\begin{split}
\int_{t_0}^{t_1} \int_M e^{-u} \brs{F}^2 \gD u_t dV dt \leq&\ \ge
\int_{t_0}^{t_1} \int_M e^u \brs{\N u_t}^2 dV dt\\
&\ + C \int_{t_0}^{t_1} \int_M \left( e^{u} \brs{\N u}^2 \brs{F}_g^4 + 
e^{2 u} \brs{\N^g F}_g^2 \brs{F}_g^2 \right) dV dt.
\end{split}
\end{gather}
We have already bounded the last term in the above inequality.  The first term
in the second line above is the one which finally requires the bound on the
Sobolev constant of $g$.  We start with an application of H\"older's inequality
and the Sobolev inequality with respect to $g_0$.
\begin{gather} \label{intest60}
\begin{split}
\int_M e^{u} \brs{\N u}^2 \brs{F}_g^4 dV \leq&\ \brs{\brs{\N u}}_{L^4}^2
\brs{\brs{e^{\frac{u}{2}} \brs{F}^2_g}}_{L^4}^2\\
\leq&\ C \brs{\brs{\N u}}_{L^2} \brs{\brs{\N u}}_{H^1} \brs{\brs{e^{\frac{u}{2}}
\brs{F}^2_g}}_{L^2} \brs{\brs{e^{\frac{u}{2}} \brs{F}^2_g}}_{H^1}\\
\leq&\ C \sup_{t_0 \leq t < t_1} \brs{\brs{u}}_{H^2} \brs{\brs{e^{\frac{u}{2}}
\brs{F}_g^2}}_{L^2} \brs{\brs{e^{\frac{u}{2}} \brs{F}_g^2}}_{H^1}.
\end{split}
\end{gather}
Now we note
\begin{align*}
\brs{\brs{e^{\frac{u}{2}} \brs{F}_g^2}}_{L^2} =&\ \left(\int_M e^u \brs{F}_g^4
dV \right)\\
=&\ \brs{\brs{ \brs{F}_g}}_{L^4(g)}^2\\
\leq&\ C_S(g) \brs{\brs{ \brs{F}_g}}_{L^2(g)} \brs{\brs{ \brs{F}_g}}_{H^1(g)}\\
\leq&\ C \left(\int_M \brs{\N^g F}_g^2 dV_g \right)^{\frac{1}{2}}\\
=&\ C \brs{\brs{e^{\frac{u}{2}} \brs{\N^g F}_g}}_{L^2}.
\end{align*}
Plugging this into (\ref{intest60}) yields
\begin{gather}
\begin{split}
\int_M e^{u} \brs{\N u}^2 \brs{F}_g^4 dV \leq&\ C \sup_{t_0 \leq t < t_1}
\brs{\brs{u}}_{H^2} \brs{\brs{e^{\frac{u}{2}} \brs{\N^g F}_g}}_{L^2}
\brs{\brs{e^{\frac{u}{2}} \brs{F}_g^2}}_{H^1}.
\end{split}
\end{gather}
Integrating this in time and arguing as in line (\ref{intest101}) yields
\begin{align*}
\int_{t_0}^{t_1} \int_M &e^{u} \brs{\N u}^2 \brs{F}_g^4 dV dt\\
\leq&\ C \ge \sup_{t_0 \leq t < t_1} \brs{\brs{u}}_{H^2}^2 + C \ge
\int_{t_0}^{t_1} \int_M \left[ e^{u} \brs{\N u}^2 \brs{F}_g^4 + e^{2 u}
\brs{\N^g F}_g^2 \brs{F}_g^2 \right] dV dt\\
\leq&\ C \ge \sup_{t_0 \leq t < t_1} \brs{\brs{u}}_{H^2}^2 + C \ge \sup_{t_0
\leq t < t_1} \int_M e^u \brs{\N^g F}^2_g dV_g + C \ge \int_{t_0}^{t_1} \int_M
e^{u} \brs{\N^g \N^g F}_g^2 dV_g.
\end{align*}
where in the last line we rearranged terms and applied (\ref{intest110}).
Thus plugging this into (\ref{intest50}), applying (\ref{intest110}) again and
plugging the result into (\ref{intest40}) gives
\begin{gather} \label{intest80}
\begin{split}
I \leq&\ C \left(t_1 - t_0 + \ge \right) I + C \brs{\brs{u(t_0)}}_{H^2}^2\\
&\ + C \ge \left( \sup_{t_0 \leq t < t_1} \int_M e^u \brs{\N^g F}_g^2 dV_g +
\int_{t_0}^{t_1} \int_M e^u \brs{\N^g \N^g F}_g^2 dV_g \right) + C(T).
\end{split}
\end{gather}
Combining this with (\ref{intest2}) and choosing $\ge$ small with respect to
universal constants gives
\begin{gather*}
\begin{split}
\sup_{t_0 \leq t < t_1} \left( \brs{\brs{u}}_{H^2}^2 + \int_M e^u \brs{\N^g
F}_g^2 dV_g \right) \leq&\ C \left( \brs{\brs{u(t_0)}}_{H^2} + \int_M e^{u}
\brs{\N^g F}_g^2 dV_g (t_0) \right)\\
&\ + C(T).
\end{split}
\end{gather*}
Thus we can cover $[0, T]$ by finitely many intervals of length $\tau$ to yield
an $H^2$ bound for $u$ and an $H^1$-type
bound for $F$ on any finite time interval.  It is easy to see that we now also
have a bound on $\brs{\brs{F}}_{H^1}$.  Now we may choose a sequence of times
$t_n \to T$ and choose
divergence-free gauges for the connections $A(t_n)$.  Our $H^1$ bound for $F$
then yields an $H^2$ bound for $A$ and so we can conclude that both $A$ and $u$
have uniform $C^{\frac{1}{2}}$ bound up to time $T$.  Using this and the form of
the evolution equations we
can apply parabolic Schauder estimates at this point to conclude $C^\infty$
convergence at $t = T$.  This completes the proof of Theorem \ref{existence} and
the existence statements of the Main Theorem.

\section{Convergence Results}
We will apply the concentration-compactness result of Chen \cite{Chen} to show
convergence of the volume-normalized flow.  Again we note that we do not require
the isoperimetric constant bound here, these statements hold for any long-time
solution of RYM flow on a surface with bounded volume.  The
statement we use is taken from \cite{struwe}.

\begin{thm} \label{compactness} (\cite{struwe} Theorem 3.1).  Let $g_n = e^{u_n}
g_0$ be a family
of smooth conformal metrics on a surface $M$ with unit volume and bounded Calabi
energy.  Then either the sequence $\{u_n \}$ is bounded in $H^2(M, g_0)$ or
there exist points $ \{x_1, \dots x_L \} \in M$ and a subsequence $\{u_n \}$
such that for any $\rho > 0$ and any $i$ we have
\begin{gather*}
\liminf_{n \to \infty} \int_{B_\rho(x_i)} \brs{K_n} dV_{g_n} \geq 2 \pi
\end{gather*}
where $K_n$ is the Gauss curvature of $g_n$.  Moreover, there holds
\begin{gather*}
2 \pi L \leq \limsup_{n \to \infty} \left( Ca(g_n) + C_0 \right)^{\frac{1}{2}} <
\infty
\end{gather*}
and either $u_n \to - \infty$ and $n \to \infty$ locally uniformly on $M \slash
\{x_1, \dots x_L \}$ or $ \{u_n \}$ is locally bounded in $H^2(M, g_0)$ away
from $\{x_1, \dots x_L \}$.
\end{thm}
 
First consider the case $\chi(M) \leq 0$.  In this case the
energy $\FF$ is bounded below.  Thus, as a consequence of Proposition
\ref{volnormenergy} we have
that
\begin{gather*}
\liminf_{t \to \infty} \ddt \FF(g(t), F(t)) = 0
\end{gather*}
Thus choose a sequence of times $\{t_n \}, t_n \to \infty$ so that
\begin{gather*}
\lim_{n \to \infty} \ddt \FF(g(t_n), F(t_n)) = 0
\end{gather*}
It is clear from (\ref{volnormenergyev}) that for this sequence we
further
have
\begin{gather} \label{convpfloc10}
\lim_{n \to \infty} \int_M \brs{\N^{g_n} F}_{g_n}^2 dV_{g_n} = 0
\end{gather}
Our goal is to show that the Calabi energy is bounded.  To do that we
expand the inner product in (\ref{volnormenergyev}).  We note that intuitively
since $\int_M \brs{\N^g F}_g^2 dV_g$ is very small, one expects that $\brs{F}^2$
is roughly parallel, so that the inner product should split.  We carry out
estimates to that effect.  First note
\begin{gather} \label{convpfloc15}
\begin{split}
- \frac{1}{2} \ddt &\FF(g(t), A(t))\\
=&\ \int_M e^{-u} \left( \gD u \right)^2 dV - 2 R_0 \int_M e^{-u} \gD u dV +
\int_M e^{-2 u} \gD u \brs{F}^2 dV\\
&\ - R_0 \int_M e^{-2 u} \brs{F}^2 dV + R_0 \int_M e^{-u} \brs{F}^2 dV\\
&\ + \frac{1}{4} \int_M e^{-3 u} \brs{F}^4 dV - \frac{1}{4} \left( \int_M e^{-u}
\brs{F}^2 \right)^2\\
\leq&\ \ge
\end{split}
\end{gather}
Also we clearly have
\begin{align*}
-2 R_0 \int_M e^{-u} \gD u dV =&\ -2 R_0 \int_M e^{-u} \brs{d u}^2 dV\\
\geq&\ 0
\end{align*}
Combining these facts, and using the uniform bound on $\int_M e^{-u} \brs{F}^2
dV$ yields
\begin{align} \label{convpfloc20}
\int_M e^{- u} \left(\gD u \right)^2 dV + \int_M e^{-2 u} \gD u \brs{F}^2 dV +
\frac{1}{4} \int_M e^{-3 u} \brs{F}^4 dV \leq C
\end{align}
We now show that the middle term here must be small which gives us the desired
bound on the Calabi energy.  In the two estimates below we will use the notation
$g$ to refer to a metric in the sequence $g(t_n)$ to simplify notation.  Fix a
small $\ge > 0$ and choose a large $n$ so that $\int_M \brs{\N^g F}_g^2 dV_g <
\ge$.  At this time we can estimate
\begin{align*}
\int_M e^{-2 u} \brs{F}^2 \gD u dV =&\ \int_M \brs{F}_g^2 \gD_g u dV_g\\
=&\ - \int_M \left< \N^g \brs{F}_g^2, \N u \right>_g dV_g\\
\leq&\ C \int_M \brs{\N^g F}_g \brs{F}_g \brs{\N u}_g dV_g\\
=&\ C \int_M \left( e^{\frac{u}{2}} \brs{\N^g F}_g \right) \left( e^{- \frac{3
u}{4}} \brs{F} \right) \left( e^{- \frac{u}{4}} \brs{\N u} \right) dV\\
\leq&\ C \left( \int_M \brs{\N^g F}_g^2 dV_g \right)^{1/2} \left( \int_M e^{-3
u} \brs{F}^4 dV \right)^{1/4} \left( \int_M e^{-u} \brs{\N u}^4 dV \right)^{1/4}
\end{align*}
Now we estimate using the Calder\'on-Zygmund inequality
\begin{align*}
\left(\int_M e^{-u} \brs{\N u}^4 dV\right)^{1/4} \leq&\ C \left(\int_M \brs{\N
u}^8 dV \right)^{1/8}\\
\leq&\ C \brs{\brs{\N u}}_{H_1^{8/5}}\\
\leq&\ C \brs{\brs{\N^2 u}}_{L^{8/5}}\\
\leq&\ C \brs{\brs{ \gD u}}_{L^{8/5}}\\
=&\ C \left(\int_M e^{\frac{4}{5} u} \left( e^{- \frac{4}{5} u} \brs{\gD
u}^{8/5} dV \right) \right)^{5/8}\\
\leq&\  C \left(\int_M e^{-u} \brs{\gD u}^2 dV \right)^{\frac{1}{2}}
\end{align*}
Plugging this into the above calculation gives
\begin{gather*}
\int_M e^{-2 u} \brs{F}^2 \gD u dV \leq C \ge \left(1 + \int_M e^{-3 u}
\brs{F}^4 dV + \int_M e^{-u} \brs{\gD u}^2 dV \right)
\end{gather*}
where $C$ is a universal constant. Thus plugging this back into
(\ref{convpfloc20}) we can
conclude that for $\ge$ small enough we have a uniform bound on
\begin{gather*}
\int_M e^{-u} \left( \gD u \right)^2 dV.
\end{gather*}
This implies that the Calabi energy, given by
\begin{gather*}
Ca(g) := \int_M \brs{K_g - \bar{K}_g}^2 dV_g = \int_M e^{-u} \brs{\gD u}^2 dV -
C_0
\end{gather*}
is bounded at these times.  Using the bound on $\brs{\brs{e^{u}}}_{L^2}$ we
have
\begin{gather*}
\begin{split}
\int_{B_{\rho}(x)} \brs{K_{g(t_n)}} dV_{g_{t_n}} \leq&\ \left(Ca(g_{t_n}) + C_0
\right)^{\frac{1}{2}} \left( \int_{B_{\rho}(x)} e^u dV_{g_0}
\right)^{\frac{1}{2}}\\
\leq&\ C \left( \int_{B_{\rho}(x)} dV_0 \right)^{\frac{1}{2}}.
\end{split}
\end{gather*}
This bound rules out the bubbling possibility of Theorem \ref{compactness}, and
so we conclude a uniform $H^2$ bound on $u$ for this sequence.  Also we have an
$H^1$ bound for $F$ as in the previous section so we can take
a
convergent subsequence,
which is in fact smoothly converging to a limit $(u_{\infty}, A_{\infty})$.  By
(\ref{convpfloc10}) we know that $F$ is covariant constant.  Thus $\brs{F}^2 =
\int_M \brs{F}_g^2 dV_g$ and so the limiting metric has constant scalar
curvature. This shows that a subsequence converges as required, but using the
nonincreasing property of $\FF$, it is clear that in fact the whole flow itself
must be converging to this metric.  

For the case $\chi(M) > 0$ we consider the gauge-fixed flow introduced in Lemma
\ref{energybound}.  Here again the energy $\FF$ is bounded below so we can argue
as above to
show that the Calabi energy for the gauge-fixed
flow is bounded for a subsequence approaching infinity.  Two terms are bounded
differently.  In particular we have
\begin{align*}
2 R_0 \int_M e^{-u} \gD u dV \leq C + \ge \int_M e^{-u} \left( \gD u \right)^2
\end{align*}
and also
\begin{align*}
R_0 \int_M e^{-2 u} \brs{F}^2 \leq C + \ge \int_M e^{-3 u} \brs{F}^4 dV.
\end{align*}
Once the Calabi energy is bounded we argue as above using Theorem
\ref{compactness} to show that a subsequence of the gauge-fixed flow converges. 
Since the
diffeomorphism parameter is defined in terms of the varying metric and we now
have uniform control and convergence of this metric, these
diffeomorphisms also converge thus the solution $(g(t), A(t))$ also converges. 
Since $F$ is parallel in the limit it follows that the limiting metric must have
constant scalar curvature.

It is clear that if one assumes a uniform upper
bound on volume for a solution to unnormalized RYM flow on $S^2$ with $[F] \neq
0$, then
using Proposition \ref{surfaceenergy}, the
arguments we have given above apply to allow us to conclude convergence to the
round metric with $F$ parallel.  This completes the proof of the main theorem.

\section{Conclusions}

Our description of RYM flow in $S^2$ is encouraging, showing
that a purely topological condition changes the qualitative behaviour of this
equation.  We have to remember however that the Ricci flow on $S^2$ always
encounters a global, type I singularity.  Indeed, an easy argument akin to lemma
\ref{volumelemma1} could show a lower volume bound for any type I singularity
when $[F^{\wedge n}] \neq 0$.  This rules out such global singularities, but
says nothing yet about local singularities, which are of course the main problem
for
Ricci flow in higher dimensions.  

We can also make the following conjecture for odd-dimensional manifolds related
to Conjecture 1:
\begin{conj} Let $(M^{2n + 1}, g)$ be a Riemannian manifold and $L \to M$
the total
space of a $U(1)$ bundle with connection $A$ satisfying $[F^{\wedge n}] \neq
0$.  Then the solution to RYM-flow exists for all time.
\end{conj}
It may be possible to use the detailed description of Ricci flows on
three-manifolds to attack this problem, in particular consider the case of $M^3
= S^2 \times S^1$.  An argument like Lemma \ref{volumelemma1} can show that the
minimal volume of an immersed $S^2$ representing the nonzero homology class can
never drop to zero.  Thus one does not expect a neckpinch singularity.  However,
at this point, even showing no local collapsing around a singularity is quite
difficult, as Perelman's proofs do not generalize in an obvious way.  Thus the
structure of singularities is still poorly understood.

Besides the lack of a bound on the isoperimetric constant, there are other
questions our main theorem leaves unanswered.  In
particular, one would like to prove that the volume of the unnormalized equation
stays bounded on $S^2$ when $[F] \neq 0$.  Also, it is likely the case that when
$[F] = 0$ the flow encounters a singularity in finite time which converges to a
round point with $F \equiv 0$.  Since resolving these questions may likely make
use of the gradient property of Ricci Yang-Mills flow, we have included the
classification of Ricci Yang-Mills solitons in dimension 2.  

It would also be interesting to know if the
flow converges exponentially.  This is likely true and would follow from a
modification of the argument in \cite{struwe}.  Finally, it would be interesting
to see if an a-priori estimate on the gradient of $u$ may be obtained similarly
to the Alexandrov reflection.  Here the presence of the Yang-Mills term in the
evolution of $u$ makes the usual proof break down.

\section{Appendix: Gradient Solitons on Surfaces}
In this section we classify all solitons on surfaces.  We note that
the proof of the corresponding result for Ricci flow using the
Kazdan-Warner identity does not work, but we can easily adapt the
proof of Chen, Lu, and Tian (\cite{CLT}, \cite{Chow}).

\begin{defn} \label{solitondef} Given $(M, g)$ a Riemannian
manifold and $L \to M$ an $S^1$ bundle over $M$ with connection $A$,
we say that $(g, A)$ is a \emph{gradient RYM soliton} if
\begin{align} \label{soliton1}
\Rc - \frac{1}{2} \eta + \N^2 f + \gl g =&\ 0\\ \label{soliton2} d^*
F =&\ \N f \hook F
\end{align}
\end{defn}

\begin{prop} \label{Schop} Ricci Yang-Mills flow is the gradient flow of the
lowest eigenvalue of the Schroedinger operator
\begin{gather}
- 4 \gD + R - \frac{1}{4} \brs{F}^2.
\end{gather}
Moreover, this eigenvalue is constant in time if and only if the solution is a
gradient Ricci Yang-Mills soliton.
\begin{proof} The proof of this proposition is adapted directly from the
corresponding proof for the Ricci flow.  It can be found in \cite{thesis}, and
was discovered independently by Andrea Young \cite{AY}.  A discussion of this
and other
gradient properties of RYM flow will appear in a subsequent paper \cite{YS}.
\end{proof}
\end{prop}

\begin{lemma} (\cite{CLT} Lemma 1) \label{sslemma1} Let $(\gS, g)$
be a two dimensional complete Riemannian manifold with non trivial
Killing vector field $X$.  If $X$ vanishes at $O \in \gS$ then
$(\gS, g)$ is rotationally symmetric
\end{lemma}

\begin{prop} \label{ssprop1} If $g$ is a gradient soliton on a closed
surface $\gS^2$ then $g$ has constant curvature and $F$ is parallel.
\begin{proof}  The gradient soliton equations on a surface
are
\begin{gather}
\begin{split}
\left(R - \frac{1}{2} \brs{F}^2 \right) g_{ij} =&\ c g_{ij} + \N_i
\N_j f\\
d^* F =&\ \N f \hook F
\end{split}
\end{gather}
for some constant $c \in \mathbb R$.  As in the case of Ricci
solitons we have that $\N f$ is a conformal vector field.  If $J$ is
the complex structure on $T \gS$ defined by counterclockwise
rotation then $J(\N f)$ is a Killing vector field, which vanishes at
some point since $\gS$ is closed.  Thus by lemma \ref{sslemma1} $g$
is rotationally symmetric. In particular we have
\begin{gather*}
g = dr^2 + \phi(r)^2 d \gt^2, \quad 0 \leq r \leq A < \infty, \quad
0 \leq \gt \leq 2 \pi
\end{gather*}
The gradient soliton equation now implies that $F$ is rotationally
symmetric also.  In particular we set $F = \psi(r) dr \wedge d \gt$
and in particular $\frac{1}{2} \brs{F}^2 = \frac{\psi^2}{\phi^{2}}$.
Now the metric component of the gradient soliton equation becomes
the pair of equations
\begin{gather}
- \frac{\phi''}{\phi} = c + \frac{\psi^2}{\phi^2} + f'', \quad -
\frac{\phi''}{\phi} = c + \frac{\psi^2}{\phi^2} + \frac{\phi'
f'}{\phi}
\end{gather}
Combining these two equations gives $f'' = \frac{\phi'}{f'}{\phi}$.
We can integrate this to give $f' = a \phi$ for a constant $a$.
Thus
\begin{gather} \label{sslemma1loc10}
- \frac{\phi''}{\phi} = c + \psi + a \phi'
\end{gather}
Next, the Yang-Mills component of the gradient soliton equation
becomes the pair of equations
\begin{gather} \label{sslemmaloc20}
\frac{\phi' \psi}{\phi} = 0, \quad \psi' = \psi f'
\end{gather}
So, multiplying (\ref{sslemma1loc10}) by $\phi \phi'$, using that
$\phi' \psi = 0$ and integrating over $[0, A]$ gives
\begin{gather*}
- c \left. \frac{(\phi')^2}{2} \right|_0^A = \left. \frac{\phi^2}{2}
\right|^{A}_0 + a \int_0^A \phi (\phi')^2 dr
\end{gather*}
Since the metric is smooth we have $\phi(0) = \phi(A) = 0$ and
$\phi'(0) = - \phi'(A) = 1$ so that $a = 0$.  Thus $f$ is constant.
By (\ref{sslemmaloc20}) we see that $\psi' = 0$, so that $F$ is
parallel.
\end{proof}
\end{prop}

\bibliographystyle{hamsplain}

\begin{thebibliography}{10}

\bibitem{Aubin} Aubin, Thierry.  \emph{Meilleures constantes dans le
th\'eor\`eme d'inclusion de Sobolev et un th\'eor\`eme de Fredholm non
lin\'eaire pour la transformation conforme de la coubure scalaire}, J. Funct.
Anal. 32 (1979), 148-174.

\bibitem{Chavel} Chavel, Isaac.  \emph{Eigenvalues in Riemannian Geometry}. 
Pure and Applied Mathematics, 115, Academic Press, Inc., Orlando, FL 1984.

\bibitem{Chen} Chen, Xiuxiong.  \emph{Calabi flow in Riemann surface revisited:
A new point of view}, Intern. Math. Res. Notices No. 6 (2001), 275-297

\bibitem{CLT} Chen, Xiuxiong; Lu, Peng; Tian, Gang. \emph{A note on
uniformization of Riemann surfaces by Ricci flow}, arXiv:math/0505163

\bibitem{Chow} Chow, Bennett; Lu, Peng; Ni, Lei. \emph{Hamilton's
Ricci Flow}, American Mathematical Society Science Press, (2006)

\bibitem{Hamilton} Hamilton, Richard. \emph{An isoperimetric estimate for the
Ricci flow on the two-sphere}.  Modern Methods in Complex Analysis (Annal of
Mathematics Studies, 137).  Princeton University Press, Princeton, NJ, 1995,
191-200.

\bibitem{Dan} Jane, Dan. Private communication.

\bibitem{Lebrun} LeBrun, C. \emph{The Einstein-Maxwell Equations, Extremal
K\"ahler Metrics and Seiberg-Witten Theory}, arXiv:0803.3734.

\bibitem{thesis} Streets, Jeffrey. \emph{Ricci Yang-Mills flow},
Ph. D. Thesis, Duke University, 2007

\bibitem{YS} Streets, Jeffrey; Young, Andrea. \emph{Gradient properties of Ricci
Yang-Mills flow}, to appear.

\bibitem{struwe} Struwe, Michael. \emph{Curvature Flows on
Surfaces}, Ann. Sc. Norm. Super. Pisa Cl. Sci. (5) 1 (2002), no. 2, 247-274.

\bibitem{wu} Wu, Hung-Hsi. \emph{The Bochner Technique in Differential
Geometry}, Mathematical Reports, 1998 Vol. 3, Harwood Academic Publishers. 

\bibitem{AY} Young, Andrea, Private communication.

\end{thebibliography}

\end{document}